\documentclass [10pt] {article}

 \usepackage [english] {babel}
 \usepackage {amssymb}
 \usepackage {amsmath}
 \usepackage {amsthm}
 \usepackage [curve] {xypic}

 \let \acuteaccent=\'
 \let \tildeaccent=\~

 \def \^{\hat}
 \def \~{\tilde}
 \def \-{\bar}
 \def \ovl{\overline}
 \renewcommand \=[1]{\bar{\bar{#1}}}
 \newcommand \Ovl[1]{\overline{\overline{#1}}}
 \def \_{\underline}
 \def \:{\colon}
 \def \%{\,{:}\,}
 \def \?{?}
 \def \le{\leq}
 \def \ge{\geq}
 \def \<{\langle}
 \def \>{\rangle}
 \def \0{\vartriangle}
 \def \${\circ}
 \def \`{\text{\rm`}}
 \def \'{\text{\rm'}}

 \def \GA{\Gamma}
 \def \Ga{\gamma}
 \def \DE{\Delta}
 \def \De{\delta}
 \def \Ep{\epsilon}
 \def \Et{\eta}
 \def \Th{\theta}
 \def \Mu{\mu}
 \def \Rh{\rho}
 \def \SI{\Sigma}
 \def \Si{\sigma}
 \def \PH{\Phi}
 \def \Ph{\phi}
 \def \Ch{\chi}

 \def \N{\boldsymbol{\mathsf N}}
 \def \Z{\boldsymbol{\mathsf Z}}
 \def \R{\boldsymbol{\mathsf R}}

 \def \d{\partial}
 \def \n{\heartsuit}

 \DeclareMathOperator {\ord} {ord}
 \DeclareMathOperator {\lk} {lk}
 \DeclareMathOperator {\id} {id}
 \DeclareMathOperator {\im} {im}
 \DeclareMathOperator {\dom} {dom}

 \newcommand* {\head} [1]
 {\subsubsection * {#1}}

 \newcommand* {\subhead} [1]
 {\addvspace\smallskipamount \noindent {\it #1\/}}

 \newenvironment* {claim} [1] []
 {\begin{trivlist}\item [\hskip\labelsep {\bf #1}] \it}
 {\end{trivlist} }

 \newenvironment* {demo} [1] []
 {\begin{trivlist}\item [\hskip\labelsep {\it #1}] }
 {\end{trivlist} }

\begin {document} \frenchspacing

 \title {\Large\bf
         Order of a homotopy invariant in the stable case}

 \author {\normalsize\rm
          Sem\"en Podkorytov}

 \date {}

 \maketitle

 \begin {abstract} \noindent
 Let
 $X$ and $Y$ be CW-complexes,
 $U$ be an abelian group, and
 $f\:[X,Y]\to U$ be a map (a homotopy invariant).
 We say that
 $f$ has {\it order\/} at most $r$ if
 the characteristic function of the $r$th Cartesian power of
 the graph of a continuous map $a\:X\to Y$ $\Z$-linearly
 determines $f([a])$.
 Suppose that
 the CW-complex $X$ is finite and
 we are in the stable case:
 $\dim X<2n-1$ and $Y$ is $(n-1)$-connected.
 We prove that then
 the order of $f$ equals
 its degree with respect to the Curtis filtration of the group
 $[X,Y]$.
 \end {abstract}


 \head {1. Introduction}

 \subhead {Order of a homotopy invariant.}
 Let $X$ and $Y$ be (topological) spaces.
 For $r\in\N$ ($=\{0,1,\dotsc\}$),
 let $E_r$ be the group of all functions $(X\times Y)^r\to\Z$.
 For a map $a\in C(X,Y)$,
 let
 $\GA_a\subset X\times Y$ be its graph and
 $I_r(a)\in E_r$ be the characteristic function of the set
 $\GA_a^r\subset(X\times Y)^r$.
 Let $D_r\subset E_r$ be the subgroup generated by the
 functions $I_r(a)$, $a\in C(X,Y)$.

 Let
 $U$ be an abelian group and
 $f\:[X,Y]\to U$ be a map.
 Define the {\it order\/} of $f$,
 $\ord f\in\^\N$ ($=\N\cup\{\infty\}$),
 to be the infimum of those $r\in\N$ for which
 there exists a homomorphism $l\:D_r\to U$ such that
 $f([a])=l(I_r(a))$ for all $a\in C(X,Y)$.
 As one easily sees,
 the existence of such $l$ for some $r$
 implies that for all greater $r$.

 \subhead {Main result.}
 Suppose that
 $X$ is a finite CW-complex,
 $Y$ is a CW-complex, and
 we are in the {\it stable case\/}:
 $\dim X\le m$,
 $Y$ is $(n-1)$-connected, and
 $m<2n-1$.
 The set $[X,Y]$ becomes an abelian group canonically.
 There is the Curtis filtration $B=(B_s)_{s=1}^\infty$,
 $[X,Y]=B_1\supset B_2\supset\dotso$, see \S~3.
 It is known \cite {Curtis} that
 $B_s=0$ for $s>2^{m-n}$.
 The degree of $f$ with respect to $B$, $\deg_Bf\in\^\N$,
 is defined, see below.

 \begin {claim} [(1.1) Theorem.]
 $\ord f=\deg_B f$.
 \end {claim}

 Example:
 if $f$ is a homomorphism,
 its order equals the greatest $s$ for which $f|B_s\ne0$.
 (If $f=0$,
 then $\ord f=0$).

 \subhead {Degree of a map between abelian groups
           with respect to a filtration.}
 Let
 $T$ and $U$ be abelian groups,
 $f\:T\to U$ be a map, and
 $P=(P_s)_{s=1}^\infty$ be a filtration of the group $T$:
 $T=P_1\supset P_2\supset\dotso$.
 Define the {\it degree\/} of $f$ with respect to $P$,
 $\deg_P f\in\^\N$,
 to be the infimum of those $r\in\N$ for which
 $$
 \sum_{e_1,\dotsc,e_k=0,1}
 (-1)^{e_1+\dotso+e_k}f(e_1t_1+\dotso+e_kt_k)=0
 $$
 whenever
 $k\in\N$,
 $t_l\in P_{s_l}$, $l=1,\dotsc,k$, and
 $s_1+\dotso+s_k>r$.


 \head {2. Preliminaries}

 \subhead {Polyhedra.}
 A {\it polyhedron\/} $L$ is
 a finite set of affine simplices in $\R^\infty$
 satisfying the ``axioms of a simplicial complex'' and
 equipped with a linear order of the vertices of each simplex
 in such a way that the order of the vertices of a simplex
 induces the order of the vertices of each of its faces.
 The {\it body\/} $|L|$ of $L$ is the union of its simplices.
 A {\it polyhedral body\/} is the body of some polyhedron.

 \subhead {Morphisms of polyhedra.}
 For polyhedra $K$ and $L$,
 a map $f\:K\to L$ is called a {\it morphism\/} if
 a vertex is sent to a vertex,
 the image of a simplex is spanned by the images of its
 vertices, and
 the non-strict order of vertices is preserved.
 A morphism $f\:K\to L$ induces a continuous map
 $|f|\:|K|\to|L|$.

 \subhead {Generation.}
 A simplex $y\in L$ generates a subpolyhedron $\=y\subset L$.
 A set $T\subset L$ generates a subpolyhedron $\-T\subset L$.

 \subhead {Small sets.}
 A set $T\subset L$ is {\it small\/} if
 there exists a simplex $y\in L$ with $\=y\supset T$;
 the least of such simplices is {\it spanned\/} by $T$.

 \subhead {The distance $\Rh_L$.}
 For $x,y\in L$,
 let $\Rh_L(x,y)\in\^\N$ be the infimum of lengths of edge
 chains connecting $x$ and $y$.
 (The orientation of edges is disregarded;
 the length of a chain is the number of its edges.)
 If $\Rh_L(x,y)<a$, $\Rh_L(y,z)<b$ ($x,y,z\in L$,
 $a,b\in\N$),
 then $\Rh_L(x,z)<a+b$.

 \subhead {Neighbourhoods $O_L$.}
 For $y\in L$ and $d\in\N$,
 put $O_L(y,d)=\{z\in L:\Rh_L(y,z)<d\}$.
 For $T\subset L$,
 let $O_L(T,d)$ be the union of the sets $O_L(y,d)$, $y\in T$.

 \subhead {Separation $\Ep_L$.}
 For $T\subset L$,
 put $\Ep_L(T)=\inf\{\Rh_L(x,y):x,y\in T,x\ne y\}\in\^\N$.

 \subhead {Subdivisions.}
 Equip the barycentric subdivision of $L$ with the following
 order:
 the greater dimension of a simplex is,
 the higher its barycentre is.
 Let $\De L$ denote the resulting polyhedron.
 Let $\Ph_L\:\De L\to L$ be the morphism taking the barycentre
 of a simplex to the highest of its vertices.
 Equip the barycentric subdivision of $L$ with the opposite
 order.
 Let $\De'L$ denote the resulting polyhedron.
 Let $\Ph'_L\:\De'L\to L$ be the morphism taking the barycentre
 of a simplex to the lowest of its vertices.
 Put
 $\DE L=\De'\De L$ and
 $\PH_L=\Ph_L\circ\Ph'_{\De L}\:\DE L\to L$.
 The map $|\PH_L|\:|L|=|\DE L|\to|L|$ is homotopic to the
 identity.
 The image of the star of each simplex of $\DE L$ under $\PH_L$
 is small.
 Thus,
 if $\Rh_{\DE L}(x,y)\le2d$ ($x,y\in\DE L$, $d\in\N$),
 then $\Rh_L(\PH_L(x),\PH_L(y))\le d$.

 \subhead {The empty simplex.}
 Put $L^\$=L\cup\{\varnothing\}$.
 Let the empty simplex generate the empty subpolyhedron:
 $\=\varnothing=\varnothing$.
 For $x,y\in L^\$$,
 we have $x\cap y\in L^\$$.

 \subhead {Completion.}
 Adding degenerate simplices to $L$,
 we get a simplicial set $\^L$.
 We have $L\subset\^L_0\cup\^L_1\cup\dotso$.
 The spaces $|L|$ and $|\^L|$ are canonically homeomorphic.
 A moprhism $f\:K\to L$ of polyhedra induces a simplicial map
 $\^f\:\^K\to\^L$.
 The correspondence $f\mapsto\^f$ is bijective.

 \subhead {Sections.}
 For a simplicial set $E$,
 let $E(L)$ be the set of simplicial maps $v\:\^L\to E$, {\it
 sections}.
 A section $v\in E(L)$ induces a map $|v|\in C(|L|,|E|)$.
 For a subpolyhedron $K\subset L$,
 we have the restriction $v|_K\in E(K)$.
 For a morphism $f\:K\to L$ of polyhedra,
 we have the composition $v\circ f\in E(K)$.
 A simplicial map $t\:D\to E$
 induces a map $t_\#\:D(L)\to E(L)$.
 For a simplicial group $G$ and a section $v\in G(L)$,
 put $\Si(v)=\{y\in L:v|_{\=y}\ne1\}$.

 \subhead {Quasisections.}
 For a set $T\subset L$ and a simplicial set $E$,
 put
 $$
 E_T=\prod_{y\in T}E(\=y).
 $$
 For $v\in E(L)$,
 put $v\|_T=(v|_{\=y})_{y\in T}\in E_T$.
 For
 a {\it quasisection\/} $w\in E_L$ and
 a morphism $f\:K\to L$ of polyhedra,
 define the composition $w\circ f\in E_K$
 by $(w\circ f)_x=w_{f(x)}\circ f'_x$, $x\in K$, where
 $f'_x\:\=x\to\Ovl{f(x)}$ are the restrictions of $f$.
 We have the map $f^\#\:E_L\to E_K$, $f^\#(w)=w\circ f$.
 For a simplicial map $t\:D\to E$ and a quasisection
 $v\in D_L$,
 we have the composition $t\circ v\in E_L$.

 \subhead {Free groups.}
 For a set $E$ with a marked element $*$,
 we have the group $FE$ given by
 the generators $\_e$, $e\in E$, and
 the relation $\_*=1$.
 The map $i\:E\to FE$, $i(e)=\_e$, is called {\it canonical}.

 \subhead {The lower central series and the abelianization.}
 For a group $G$,
 let $(\Ga_sG)_{s=1}^\infty$ be its lower central series.
 Put $G^+=G/\Ga_2G$.

 \subhead {Free abelian groups.}
 For a set $E$,
 we have the abelian group $\<E\>$ with the base
 $(\`e\')_{e\in E}$.
 The map $j\:E\to\<E\>$, $j(e)=\`e\'$, is called {\it
 canonical}.
 Let $\<E\>_\0$ be the kernel of the homomorphism $\<E\>\to\Z$,
 $\`e\'\mapsto1$.
 A map $t\:D\to E$ induces a homomorphism
 $\<t\>\:\<D\>\to\<E\>$.

 Let
 $L$ be a polyhedron,
 $E$ be a simplicial set, and
 $V\in\<E(L)\>$ be an element (an {\it ensemble\/}).
 Let $|V|\in\<C(|L|,|E|)\>$ denote the image of $V$ under the
 homomorphism induced by the map $|?|\:E(L)\to C(|L|,|E|)$.
 For a subpolyhedron $K\subset L$,
 the ensemble $V|_K\in\<E(K)\>$ is defined similarly;
 for a set $T\subset L$,
 we have the element $V\|_T\in\<E_T\>$.
 For spaces $X$ and $Y$ and an ensemble $A\in\<C(X,Y)\>$,
 we have the element $[A]\in\<[X,Y]\>$.
 For a set $Z\subset X$,
 we have the ensemble $A|_Z\in\<C(Z,Y)\>$.

 For a simplicial group $G$ and an ensemble $V\in\<G(L)\>$,
 $$
 V=\sum_{v\in G(L)}m_v\`v\'
 $$
 ($m_v\in\Z$),
 put
 $$
 \SI(V)=\bigcup_{v\in G(L)\%m_v\ne0}\Si(v).
 $$

 \subhead {Group rings.}
 For a group $G$,
 $\<G\>$ is the group ring,
 $\<G\>_\0$ is its (two-sided) ideal.
 For $s\in\N_+$ ($=\N\setminus\{0\}$),
 the ideal $\<G\>_\0^s$ is additively generated by all elements
 of the form $(\`g_1\'-1)\dotso(\`g_s\'-1)$,
 $g_1,\dotsc,g_s\in G$.

 \subhead {Simplicial application.}
 Natural constructions can be applied to simplicial objects
 dimension-wise.
 For a pointed simplicial set $E$,
 we have
 the simplicial group $FE$ and
 the canonical simplicial map $i\:E\to FE$.
 The map $i$ is a model of the canonical map of a pointed space
 to the loop space of its suspension ({\it Milnor's model}, see
 \cite {Milnor}).
 For a simplicial group $G$,
 we have
 the simplicial abelian group $G^+$,
 the simplicial ring $\<G\>$,
 the canonical simplicial map $j\:G\to\<G\>$, and
 the simplicial subgroups
 $\Ga_sG\subset G$, $s\in\N_+$, and
 $\<G\>_\0^s\subset\<G\>$, $s\in\N$.

 \subhead {Simplicial trifles.}
 A simplicial map between pointed simplicial sets is called
 {\it bound\/}
 if it preserves the pointing.
 A simplicial abelian group $D$ is called {\it free\/}
 if the abelian groups $D_n$, $n\in\N$, are free.
 For a simplicial set $E$,
 let $E_{(m)}\subset E$ ($m\in\N$) denote its $m$-skeleton.

 \subhead {Fusion.}
 Let
 $L$ be a polyhedron and
 $G$ be a simplicial group.
 Let $j\:G\to\<G\>$ be the canonical map.
 The ring homomorphism $J\:\<G(L)\>\to\<G\>(L)$,
 $J(\`v\')=j\circ v$, is called {\it fusion}.


 \head {3. The Curtis filtration in the stable case}

 Let $X$ and $Y$ be CW-complexes.
 Suppose that
 $\dim X\le m$,
 $Y$ is $(n-1)$-connected, and
 $m<2n-1$.
 We shall construct a filtration $B=(B_s)_{s=1}^\infty$ of the
 abelian group $[X,Y]$, $[X,Y]=B_1\supset B_2\supset\dotso$,
 the {\it Curtis filtration}.
 There are
 a simplicial set $E$ and
 a homotopy equivalence $k\:Y\to|E|$.
 Let us point $E$.
 We have the simplicial group $G=FE$.
 By the Freudenthal theorem,
 the canonical simplicial map $i\:E\to G$ is
 $(2n-1)$-connected.
 The map $h=|i|\circ k\:Y\to|G|$ is also $(2n-1)$-connected.
 Let $j_s\:\Ga_sG\to G$, $s\in\N_+$, be the inclusions.
 For $s\in\N_+$,
 we have the chain of groups and homomorphisms
 $$
 \xymatrix {
 [X,Y]
 \ar[r]^-{h_*} &
 [X,|G|] &
 [X,|\Ga_sG|].
 \ar[l]_-{|j_s|_*}
 }
 $$
 Since $m<2n-1$,
 $h_*$ is an isomorphism.
 Put $B_s=h_*^{-1}(\im|j_s|_*)$.
 (The result does not depend on the choice of $E$ etc.)


 \head {4. A claim on Lie rings}

 Here $U$ denotes the universal enveloping ring functor.

 \begin {claim} [(4.1)]
 Let
 $L$ and $M$ be Lie rings, free as abelian groups, and
 $k\:L\to M$ be an injective homomorphism.
 Then the homomorphism $Uk\:UL\to UM$ is injective.
 \end {claim}

 \begin {demo}
 This follows easily from the
 Poincar\acuteaccent{e}--Birkhoff--Witt theorem.
 \qed
 \end {demo}


 \head {5. A claim on group rings}

 Let
 $V$ and $W$ be groups and
 $t\:V\to W$ be a homomorphism.
 We have the ring homomorphism $\<t\>\:\<V\>\to\<W\>$.
 For $s\in\N$,
 let $I_s\subset\<V\>$ be the subgroup generated by all
 elements of the form $(\`v_1\'-1)\dotso(\`v_k\'-1)$, where
 $k\in\N$,
 $v_l\in t^{-1}(\Ga_{s_l}W)$, and
 $s_1+\dotso+s_k\ge s$.
 It is easy to see that
 $I_s$ are ideals,
 $I_s\supset I_{s+1}$, and
 $I_sI_t\subset I_{s+t}$.

 \begin {claim} [(5.1)]
 Suppose that $W$ is a product of a finite number of free
 groups.
 Then $\<t\>^{-1}(\<W\>_\0^s)=I_s$, $s\in\N$.
 \end {claim}

 \begin {demo} [Proof.]
 If $w\in\Ga_sW$,
 then $\`w\'-1\in\<W\>_\0^s$
 (this holds for arbitrary $W$ \cite [III.1.3] {Passi}).
 This yields the inclusion $\<t\>^{-1}(\<W\>_\0^s)\supset I_s$.

 We have the graded rings
 $P$, $P_s=I_s/I_{s+1}$, and
 $Q$, $Q_s=\<W\>_\0^s/\<W\>_\0^{s+1}$.
 Since $\<t\>(I_s)\subset\<W\>_\0^s$,
 the homomorphism $\<t\>$ induces a graded ring homomorphism
 $l\:P\to Q$.
 We shall show that $l$ is injective.
 Then induction on $s$ with application of the 5-lemma shows
 that the induced homomorphism $\<V\>/I_s\to\<W\>/\<W\>_\0^s$
 is injective,
 which is the desired equality.

 We have the graded Lie rings
 $L$, $L_s=t^{-1}(\Ga_sW)/t^{-1}(\Ga_{s+1}W)$, and
 $M$, $M_s=\Ga_sW/\Ga_{s+1}W$
 (the product is induced by the group commutator,
 see \cite [VIII.2] {Passi}).
 The homomorphism $t$ induces a graded Lie ring homomorphism
 $k\:L\to M$,
 which is obviously injective.

 We have the commutative diagram
 $$
 \xymatrix {
 L
 \ar[r]^-{k}
 \ar[d]_-{f} &
 M
 \ar[d]^-{g} \\
 P
 \ar[r]^-{l} &
 Q,
 }
 $$
 where $f$ and $g$ are the representations with the components
 $f_s\:L_s\to P_s$, $f_s(v)=\`v\'-1$, $v\in t^{-1}(\Ga_sW)$,
 and $g_s\:M_s\to Q_s$, $g_s(w)=\`w\'-1$, $w\in\Ga_sW$.
 Extending the representations $f$ and $g$ to homomorphisms of
 the universal enveloping rings,
 we get the commutative diagram
 $$
 \xymatrix {
 UL
 \ar[r]^-{Uk}
 \ar[d]_-{\~f} &
 UM
 \ar[d]^-{\~g} \\
 P
 \ar[r]^-{l} &
 Q.
 }
 $$
 By Magnus' method, one easily shows that
 $\~g$ is an isomorphism, and
 $M$ is free as an abelian group
 (cf. \cite [VIII.6] {Passi}).
 By (4.1), 
 the homomorphism $Uk$ is injective.
 The ring $P$ is generated by elements of the form
 $\`v\'-1\in P_s$, where $s\in\N_+$, $v\in t^{-1}(\Ga_sW)$.
 They belong to the image of the representation $f$ and,
 consequently, of the homomorphism $\~f$,
 which is thus surjective.
 Therefore,
 the homomorphism $l$ is injective
 (and $\~f$ is an isomorphism.)
 \qed
 \end {demo}


 \head {6. Some ideals of the group ring of a product of
        groups}

 Let $(G_i)_{i\in I}$ be a finite collection of groups.
 For $J\subset I$,
 put
 $$
 G_J=\prod_{i\in J}G_i,
 $$
 and let $p_J\:G_I\to G_J$ be the projection homomorphism.
 We have the ring homomorphisms $\<p_J\>\:\<G_I\>\to\<G_J\>$.

 \begin {claim} [(6.1)]
 For $s\in\N$,
 we have
 $$
 \bigcap_{\#J<s}\ker\<p_J\>\subset\<G_I\>_\0^s.
 $$
 \end {claim}

 \begin {demo} [Proof.]
 We have
 $$
 \<G_I\>=\bigotimes_{i\in I}\<G_i\>.
 $$
 Since $\<G_i\>=\<G_i\>_\0\oplus\<1\>$,
 $$
 \<G_I\>=\bigoplus_{J\subset I}S(J), \qquad
 S(J)=\bigotimes_{i\in I}T_i(J),
 $$
 where the subgroup $T_i\subset\<G_i\>$ is:
 $\<G_i\>_\0$ if $i\in J$, and
 $\<1\>$ otherwise.
 Obviously,
 $\<p_J\>|S(J')$ is:
 a monomorphism if $J'\subset J$, and
 zero otherwise.
 Therefore,
 $$
 \bigcap_{\#J<s}\ker\<p_J\>=\bigoplus_{\#J\ge s}S(J).
 $$
 Now it suffices to note that $S(J)\subset\<G_I\>_\0^{\#J}$.
 \qed
 \end {demo}


 \head {7. The functions $\Et$ and $\Th$}

 Let
 $L$ be a polyhedron and
 $G$ be a simplicial group.
 We have the homomorphism $\?\|_L\:G(L)\to G_L$.
 For $V\in\<G(L)\>$,
 put $\Et(V)=\sup\{s\in\N:V\|_L\in\<G_L\>_\0^s\}\in\^\N$.
 For $s\in\N$,
 we have the subgroup $I_s\subset\<G(L)\>$
 generated by all elements of the form
 $(\`v_1\'-1)\dotso(\`v_k\'-1)$, where
 $k\in\N$,
 $v_l\in(\Ga_{s_l}G)(L)\subset G(L)$, and
 $s_1+\dotso+s_k\ge s$.
 (It is an ideal.)

 \begin {claim} [(7.1)]
 Suppose that the groups $G_n$, $n\in\N$, are free.
 Then $\{V\in\<G(L)\>:\Et(V)\ge s\}=I_s$, $s\in\N$.
 \end {claim}

 \begin {demo}
 This follows from (5.1).
 \qed
 \end {demo}

 For a simplicial set $E$ and an ensemble $V\in\<E(L)\>$,
 put $\Th(V)=\inf\{\#T:T\subset L,V\|_T\ne0\}\in\^\N$.

 \begin {claim} [(7.2)]
 For $V\in\<G(L)\>$,
 we have $\Th(V)\le\Et(V)$.
 \end {claim}

 \begin {demo}
 This follows from (6.1).
 \qed
 \end {demo}


 \head {8. Product of affine functions}

 \begin {claim} [(8.1)]
 Let
 $V$ be a group,
 $H$ be a ring, and
 $a_1,\dotsc,a_r\:V\to H$ be homomorphisms
 (to the additive group;
 $r\in\N$).
 We have the additive homomorphism $Q\:\<V\>\to H$,
 $$
 Q(\`v\')=\prod_{s=1}^r(1+a_s(v)).
 $$
 Then $Q|\<V\>_\0^{r+1}=0$.
 \end {claim}

 \begin {demo}
 This follows from \cite [V.2.1] {Passi}.
 \qed
 \end {demo}


 \head {9. Strict and $r$-strict homomorphisms}

 Let $V$ and $W$ be groups.
 An additive homomorphism $h\:\<V\>\to\<W\>$ is called
 {\it strict\/} if $h(\<V\>_\0^s)\subset\<W\>_\0^s$ for all
 $s\in\N$ and
 {\it $r$-strict\/} ($r\in\N$) if this holds for $s\le r$.

 \begin {claim} [(9.1)]
 Let $t\:V\to W$ be a homomorphism.
 Then the homomorphism $\<t\>\:\<V\>\to\<W\>$ is strict.
 \qed
 \end {claim}

 \begin {claim} [(9.2)]
 Let $f,g\:\<V\>\to\<W\>$ be $r$-strict ($r\in\N$)
 homomorphisms.
 Then the homomorphism $h\:\<V\>\to\<W\>$,
 $h(\`v\')=f(\`v\')g(\`v\')$, is $r$-strict.
 \end {claim}

 \begin {demo} [Proof.]
 Take
 $s\in\N_+$, $s\le r$, and
 $v_1,\dotsc,v_s\in V$.
 Put $x_t=\`v_t\'-1\in\<V\>_\0$.
 Let us show that $h(x_1\dotso x_s)\in\<W\>_\0^s$.
 We have
 \begin {multline*}
 (-1)^sh(x_1\dotso x_s)=
 \sum_{e_1,\dotsc,e_s=0,1}
 (-1)^{e_1+\dotso+e_s}
 h(\`v_1^{e_1}\dotso v_s^{e_s}\')= \\ =
 \sum_{e_1,\dotsc,e_s=0,1}
 (-1)^{e_1+\dotso+e_s}
 f(\`v_1^{e_1}\dotso v_s^{e_s}\')
 g(\`v_1^{e_1}\dotso v_s^{e_s}\') = \\ =
 \sum_{e_1,\dotsc,e_s=0,1}
 (-1)^{e_1+\dotso+e_s}
 f\bigl(\prod_{t=1}^s(1+e_tx_t)\bigr)
 g\bigl(\prod_{t=1}^s(1+e_tx_t)\bigr) = \\ =
 \sum_{e_1,\dotsc,e_s=0,1}
 (-1)^{e_1+\dotso+e_s}
 \bigl(
 \sum_{a_1,\dotsc,a_s=0,1}
 e_1^{a_1}\dotso e_s^{a_s}
 f(x_1^{a_1}\dotso x_s^{a_s})
 \bigr)\cdot \\ \hfill \cdot
 \bigl(
 \sum_{b_1,\dotsc,b_s=0,1}
 e_1^{b_1}\dotso e_s^{b_s}
 g(x_1^{b_1}\dotso x_s^{b_s})
 \bigr)= \\ =
 \sum_{a_1,b_1,\dotsc,a_s,b_s=0,1}
 \bigl(
 \sum_{e_1,\dotsc,e_s=0,1}
 (-1)^{e_1+\dotso+e_s}
 e_1^{a_1+b_1}\dotso e_s^{a_s+b_s}
 \bigr)
 f(x_1^{a_1}\dotso x_s^{a_s}) \cdot \\ \cdot
 g(x_1^{b_1}\dotso x_s^{b_s}).
 \end {multline*}
 Fix $a_1,b_1,\dotsc,a_s,b_s$.
 We show that the corresponding summand of the outer sum
 belongs to $\<W\>_\0^s$.
 Put
 $a=a_1+\dotso+a_s$,
 $b=b_1+\dotso+b_s$.
 Since $a,b\le s\le r$ and the homomorphisms $f$ and $g$ are
 $r$-strict,
 we have
 $$
 f(x_1^{a_1}\dotso x_s^{a_s})
 g(x_1^{b_1}\dotso x_s^{b_s})\in\<W\>_\0^{a+b}.
 $$
 If $a+b\ge s$,
 this suffices.
 Otherwise,
 there is $t$ such that $a_t=b_t=0$.
 Then
 the quantity $e_1^{a_1+b_1}\dotso e_s^{a_s+b_s}$ does not
 depend on $e_t$, and
 thus the inner sum equals zero.
 \qed
 \end {demo}


 \head {10. Group ring of a free group}

 Let $E$ be a pointed set.
 Put $G=FE$.
 Let $i\:E\to G$ be the canonical map.
 For $s\in\N$,
 we have
 the ponted set $E^{\wedge s}=E\wedge\dotso\wedge E$
 ($E^{\wedge0}$ is the 0-sphere) and
 the homomorphism $k_s\:\<E^{\wedge s}\>_\0\to\<G\>_\0^s$,
 $$
 k_s(\`(e_1,\dotsc,e_s)\'-\`*\')=\prod_{t=1}^s(\`\_{e_t}\'-1),
 $$
 where $*\in E^{\wedge s}$ is the marked element.
 By \cite [VIII.6.2] {Passi},
 the composition
 $$
 \xymatrix {
 \<E^{\wedge s}\>_\0
 \ar[r]^-{k_s} &
 \<G\>_\0^s
 \ar[rr]^-{\text{\rm projection}} & &
 \<G\>_\0^s/\<G\>_\0^{s+1}
 }
 $$
 is an isomorphism.
 Therefore,
 $\<G\>_\0^s=D^s\oplus\<G\>_\0^{s+1}$, where
 $D^s\cong\<E^{\wedge s}\>_\0$.


 \head {11. Lift of a simplicial homomorphism}

 \begin {claim} [(11.1)]
 Consider the diagram
 $$
 \xymatrix {
 &
 Q
 \ar[d]^-{f} \\
 D
 \ar[r]^-{s} &
 P
 }
 $$
 of simplicial abelian groups and homomorphisms.
 Suppose that
 $D$ is free and $m$-connected ($m\in\N$) and
 $f$ is surjective.
 Then there exists a simplicial homomorphism $t\:D\to Q$ such
 that $f\circ t|D_{(m)}=s|D_{(m)}$.
 \end {claim}

 \begin {demo} [Proof.]
 Let $\n$ denote the normalization functor.
 The complex $D^\n$ is free.
 Thus
 $D^\n=C^0\oplus C^1\oplus\dotso$, where
 $C^n$ is a free complex with $C^n_i=0$ for $i\ne n,n+1$ and
 the differential $\d\:C^n_{n+1}\to C^n_n$ injective.
 The complex $D^\n$ is $m$-connected.
 Thus,
 for $n\le m$,
 the differential $\d\:C^n_{n+1}\to C^n_n$ is an isomorphism.
 The morphism $f^\n\:Q^\n\to P^\n$ is surjective.
 Thus,
 for $n\le m$,
 there is a morphism $g^n\:C^n\to Q^\n$ such that
 $f^\n\circ g^n=s^\n|C^n$.
 We have the morphism $h\:D^\n\to Q^\n$ with $h|C^n$ equal to:
 $g^n$ if $n\le m$, and
 zero otherwise.
 Obviously, $(f^\n\circ h)_n=s^\n_n$ for $n\le m$.
 The Dold--Kan correspondence yields the simplicial
 homomorphism $t\:D\to Q$ with $t^\n=h$.
 It has the desired property.
 \qed
 \end {demo}


 \head {12. The function $\Mu_L$}

 Let $L$ be a polyhedron.
 For $x\in L^\$$,
 put $\Mu_L(x)=1-\Ch(\lk_Lx)$
 ($\Ch$ is the Euler characteristic;
 $\lk$ is the link;
 convention: $\lk_L\varnothing=L$).

 \begin {claim} [(12.1)]
 For $y,z\in L^\$$,
 we have
 $$
 \sum_{x\in L^\$\%x\cap y=z}\Mu_L(x)=
 \begin {cases}
 1 &
 \text {if $y=z$}, \\
 0 &
 \text {otherwise}.
 \end {cases}
 $$
 \end {claim}

 \begin {demo} [Proof.]
 For $t\in L^\$$,
 we have
 $$
 \sum_{x\in L^\$\%x\subset t,\,x\cap y=z}(-1)^{\dim x}=
 \begin {cases}
 (-1)^{\dim z} &
 \text {if $z\subset t\subset y$}, \\
 0 &
 \text {otherwise}
 \end {cases}
 $$
 (convention: $\dim\varnothing=-1$).
 For $x\in L^\$$,
 we have
 $$
 \Ch(\lk_Lx)=
 \sum_{t\in L^\$\%x\varsubsetneq t}(-1)^{\dim t-\dim x-1},
 $$
 and thus
 $$
 \Mu_L(x)=
 \sum_{t\in L^\$\%x\subset t}(-1)^{\dim t-\dim x}.
 $$
 We have
 \begin {multline*}
 \sum_{x\in L^\$\%x\cap y=z}
 \Mu_L(x)=
 \sum_{x,t\in L^\$\%x\subset t,\,x\cap y=z}
 (-1)^{\dim t-\dim x}= \\ =
 \sum_{t\in L^\$}
 (-1)^{\dim t}
 \sum_{x\in L^\$\%x\subset t,\,x\cap y=z}
 (-1)^{\dim x}= \\ =
 \sum_{t\in L^\$\%z\subset t\subset y}
 (-1)^{\dim t+\dim z}=
 \begin {cases}
 1 &
 \text {if $y=z$}, \\
 0 &
 \text {otherwise}.
 \end {cases}
 \end {multline*}
 \qed
 \end {demo}


 \head {13. Dummy of a simplicial group}

 \subhead {A model of the path fibration.}
 Let $B$ be the cosimplicial simplicial pointed set where
 $B^n_m$ is the set of non-strictly increasing partial maps
 $b\:[m]\dashrightarrow[n]$ (we have $\dom b\subset[m]$)
 with the marked element $o^n_m$, $\dom o^n_m=\varnothing$, and
 the structure maps are obvious.
 For $n\in\N$,
 we have the pointed simplicial set $B^n$.

 Let $G$ be a simplicial group.
 Let $\~G$, the {\it dummy}, be the simplicial group where
 $\~G_n$ is the group of bound simplicial maps $B^n\to G$ and
 the structure homomorphisms are induced by the cosimplicial
 structure.

 \begin {claim} [(13.1)]
 The space $|\~G|$ is contractible.
 \end {claim}

 \begin {demo} [Proof.]
 Let $I$ be the simplicial set that is the standard 1-simplex:
 $I_n$ is the set of non-strictly increasing maps $s\:[n]\to[1]$.
 The collection of maps $I_n\times B^n_m\to B^n_m$,
 $(s,b)\mapsto b|(s\circ b)^{-1}(1)$, $m,n\in\N$, induces
 a contracting homotopy $I\times\~G\to\~G$.
 \qed
 \end {demo}

 Evaluation at the elements $i_n\in B^n_n$, $i_n=\id\:[n]\to[n]$,
 yields the simplicial homomorphism $p\:\~G\to G$, the {\it
 projection}.

 \begin {claim} [(13.2)]
 Suppose that $G_0=1$.
 Then $p$ is surjective.
 \end {claim}

 \begin {demo} [Proof.]
 Take an element $g\in G_n$ ($n\in\N$).
 We seek an element $\~g\in\~G_n$ with $p_n(\~g)=g$, that is,
 a bound simplicial map $\~g\:B^n\to G$ with $\~g_n(i_n)=g$.
 Let $V\subset B^n$ be the simplicial subset generated by the
 elements $i_n$ and $l_n\in B^n_1$, $\dom l_n=\{0\}$, $l_n(0)=0$.
 It is the wedge of the standard $n$-simplex and 1-simplex.
 We have the simplicial map $f\:V\to G$, $f_n(i_n)=g$,
 $f_1(l_n)=1$.
 Since
 $V$ is contractible and
 $G$ is a Kan set,
 $f$ extends to $B^n$,
 which yields the desired $\~g$.
 \qed
 \end {demo}

 \subhead {Extension of sections.}
 Let $L$ be a polyhedron.
 Take simplices $x,y\in L$ of dimensions $r$, $s$, respectively.
 Let $i\:[r]\to L$ and $j\:[s]\to L$ be the increasing
 enumerations of their vertices.
 We have the partial map
 $t=i^{-1}\circ j\:[s]\dashrightarrow[r]$.
 For a bound simplicial map $\~g\:B^r\to G$,
 let $e_{xy}(\~g)\:B^s\to G$ be the bound simplicial map such
 that
 $e_{xy}(\~g)_m(b)=\~g_m(t\circ b)$
 for $b\:[m]\dashrightarrow[s]$ ($m\in\N$).
 Thus we have the homomorphism $e_{xy}\:\~G_r\to\~G_s$.

 For $x\in L$, $\dim x=r$,
 let the homomorphism $E_x\:\~G(\=x)\to\~G(L)$ be given by
 $E_x(v)_s(y)=e_{xy}(v_r(x))$ ($y\in L$, $\dim y=s$).
 Extend this construction to the case $x\in L^\$$:
 put $E_\varnothing(1)=1$
 (we have $\~G(\=\varnothing)=1$).

 \begin {claim} [(13.3)]
 For $x\in L^\$$ and $v\in\~G(\=x)$,
 we have
 \begin {itemize}
 \item [(a)]
 $E_x(v)|_{\=x}=v$;
 \item [(b)]
 $E_x(v)|_{\=y}=E_{x\cap y}(v|_{\=x\cap\=y})|_{\=y}$
 ($y\in L^\$$);
 \item [(c)]
 $\Si(E_x(v))\subset O_L(x,1)$
 if $x\neq\varnothing$.
 \qed
 \end {itemize}
 \end {claim}

 \subhead {Realization.}
 Let
 $\~J\:\<\~G(L)\>\to\<\~G\>(L)$ and
 $\~J_x\:\<\~G(\=x)\>\to\<\~G\>(\=x)$, $x\in L$,
 be fusions.
 Obviously, $\~J_x$ are isomorphisms.
 We have the additive homomorphism,
 the {\it realization},
 $R\:\<\~G\>(L)\to\<\~G(L)\>$,
 $$
 R(w)=
 \sum_{x\in L}\Mu_L(x)(\<E_x\>\circ\~J_x^{-1})(w|_{\=x}).
 $$
 We have $R(\<\~G\>_\0(L))\subset\<\~G(L)\>_\0$.

 \begin {claim} [(13.4)]
 For $w\in\<\~G\>_\0(L)$,
 we have $\~J(R(w))=w$.
 \end {claim}

 \begin {demo} [Proof.]
 For $z\in L^\$$,
 we have the homomorphism $H_z\:\<\~G\>(\=z)\to\<\~G(L)\>$ with
 $H_z=\<E_z\>\circ\~J_z^{-1}$, $z\neq\varnothing$, and
 $H_\varnothing=0$.
 It follows from (13.3~b) that
 for $x,y\in L^\$$ and $u\in\<\~G\>_\0(\=x)$,
 we have $H_x(u)|_{\=y}=H_{x\cap y}(u|_{\=x\cap\=y})|_{\=y}$.
 For $y\in L$,
 we have
 \begin {multline*}
 \~J(R(w))|_{\=y}=
 \~J_y(R(w)|_{\=y})= \\ =
 \sum_{x\in L^\$}
 \Mu_L(x)
 \~J_y(H_x(w|_{\=x})|_{\=y})=
 \sum_{x\in L^\$}
 \Mu_L(x)
 \~J_y(H_{x\cap y}(w|_{\=x\cap\=y})|_{\=y})= \\ =
 \sum_{z\in L^\$}
 \bigl(
 \sum_{x\in L^\$\%x\cap y=z}
 \Mu_L(x)
 \bigr)
 \~J_y(H_z(w|_{\=z})|_{\=y})
 \;\overset{\text{by (12.1)}}=\;
 \~J_y(H_y(w|_{\=y})|_{\=y})= \\ =
 \~J_y(\<E_y\>(\~J_y^{-1}(w|_{\=y}))|_{\=y})
 \;\overset{\text{by (13.3~a)}}=\;
 \~J_y(\~J_y^{-1}(w|_{\=y}))=
 w|_{\=y}.
 \qed
 \end {multline*}
 \end {demo}

 \begin {claim} [(13.5)]
 For $w\in\<\~G\>(L)$,
 we have $\SI(R(w))\subset O_L(\Si(w),1)$.
 \end {claim}

 \begin {demo}
 This follows from (13.3~c).
 \qed
 \end {demo}

 \begin {claim} [(13.6)]
 We have $R(\<\~G\>_\0^s(L))\subset\<\~G(L)\>_\0^s$, $s\in\N$.
 \end {claim}

 \begin {demo} [Proof.]
 For $w\in\<\~G\>_\0^s(L)$ and $x\in L$,
 we have
 $w|_{\=x}\in\<\~G\>_\0^s(\=x)$,
 $\~J_x^{-1}(w|_{\=x})\in\<\~G(\=x)\>_\0^s$, and,
 by (9.1),
 $\<E_x\>(\~J_x^{-1}(w|_{\=x}))\in\<\~G(L)\>_\0^s$.
 Summing over $x\in L$,
 we get $R(w)\in\<\~G(L)\>_\0^s$.
 \qed
 \end {demo}


 \head {14. Partitions}

 Let
 $L$ be a polyhedron and
 $D$ be a simplicial abelian group.
 A collection $(h_z\:D(\=z)\to D(L))_{z\in L}$ of homomorphisms
 is called a {\it partition\/} if
 for $w\in D(L)$,
 we have
 $$
 \sum_{z\in L}h_z(w|_{\=z})=w
 $$
 and
 $\Si(h_z(w))\subset O_L(z,1)$ for all $z\in L$.

 \begin {claim} [(14.1)]
 Suppose that
 $\dim L\le m$ ($m\in\N$) and
 $D$ is free and $m$-connected.
 Then there exists a partition
 $(h_z\:D(\=z)\to D(L))_{z\in L}$.
 \end {claim}

 \begin {demo} [Proof.]
 We shall use the Dold--Kan correspondence.
 There is a decomposition $D=D^0\oplus D^1\oplus\dotso$, where
 $D^n$ is a simplicial abelian group such that
 its normalization $C^n$ is concentrated in dimensions $n$ and
 $n+1$ and
 the differential $\d\:C^n_{n+1}\to C^n_n$ is injective
 (cf. proof of (11.1)).
 It suffices to construct a partition
 $(h^n_z\:D^n(\=z)\to D^n(L))_{z\in L}$ for each $n$.
 Take $n\le m$.
 Then $\d\:C^n_{n+1}\to C^n_n$ is an isomorphism,
 since $D$ is $m$-connected.
 Thus a section on a polyhedron with values in $D^n$ is the same
 as an $n$-cochain on it with coefficients in $C^n_n$.
 Let $h^n_z$ be:
 the extension of a cochain by zero if $\dim z=n$, and
 zero otherwise.
 Take $n>m$.
 Then
 $D^n(L)=0$
 since $\dim L\le m$.
 Thus there is the zero partition.
 \qed
 \end {demo}


 \head {15. Modification of an ensemble of sections}

 Fix numbers $b_1,\dotsc,b_5,c\in\N$ such that
 each is sufficiently great with respect to the previous,
 namely:
 $b_1\ge2$,
 $b_2\ge b_1+2$,
 $b_3\ge2b_2$,
 $b_4\ge2b_1+b_3$,
 $b_5\ge2b_2+b_4$,
 $2^{c-1}\ge2b_5+1$.

 \subhead {The morphism $e\:L\to K$.}
 Let $K$ be a polyhedron with $\dim K\le m$ ($m\in\N$).
 Put
 $L=\DE^cK$ and
 $e=\PH_K\circ\dotso\circ\PH_{\DE^{c-1}K}\:L\to K$.
 For $z\in L$,
 the set $e(O_L(z,b_5))\subset K$ is small
 (this follows from
 the properties of the operation $\DE$ and
 the inequaity $2^{c-1}\ge2b_5+1$).

 \subhead {The morphisms $e_z$.}
 Take a simplex $z\in L$.
 Since $b_2\le b_5$,
 the set $e(O_L(z,b_2))$ is small.
 It spans a simplex $x\in K$.
 Let $u\in K$ be the highest vertex of $x$.
 We shall construct a morphism $e_z\:L\to K$ with the following
 properties:
 \begin {itemize}
 \item [(1)]
 $e_z(O_L(z,b_1))=\{u\}$;
 \item [(2)]
 $e_z(O_L(z,b_2))\subset\=x$;
 \item [(3)]
 $e_z$ agrees with $e$ outside $O_L(z,b_2)$.
 \end {itemize}
 Put $L_1=\De\DE^{c-1}K$.
 We have $L=\De'L_1$.
 Let $B_1\subset L_1$ be the subpolyhedron
 generated by the simplices
 whose centres
 (which are vertices of $L$)
 belong to $O_L(z,b_1+1)$.
 Put $B=\De'B_1$.
 We have $B\subset L$ (a subpolyhedron).
 We have
 $O_L(z,b_1)\subset B$ and
 (since $b_2\ge b_1+2$)
 $O_L(B,1)\subset O_L(z,b_2)$.
 The polyhedron $L$ has no edges outcoming from $B$.
 Let $e_z$ take a vertex $t\in L$ to:
 $u$ if $t\in B$, and
 $e(t)$ otherwise.
 One easily checks that $e_z$
 is well-defined and
 has the desired properties.

 \subhead {The morphisms $e_Z$.}
 Take a set $Z\subset L$ with $\Ep_L(Z)\ge b_3$.
 Define a morphism $e_Z\:L\to K$ by the following conditions:
 \begin {itemize}
 \item [(1)]
 for $z\in Z$,
 the morphisms $e_Z$ and $e_z$ agree on $O_L(z,b_2)$;
 \item [(2)]
 the morphisms $e_Z$ and $e$ agree outside $O_L(Z,b_2)$.
 \end {itemize}
 Since $b_3\ge2b_2$,
 $e_Z$ is well-defined.

 \subhead {The simplicial groups $G$ and $D$.}
 Let $E$ be an $(n-1)$-connected ($n\in\N$) simplicial set with
 a single vertex.
 Suppose that $m\le2n-1$.
 Put $G=FE$.
 Let
 $i\:E\to G$ and $j\:G\to\<G\>$ be the canonical simplicial maps and
 $q\:G\to G^+$ be the simplicial homomorphism that is the projection.
 We shall need a decomposition
 $\<G\>\cong\<1\>\oplus G^+\oplus D$ (cf. \S~10) and some
 related simplicial homomorphisms.
 Let $d\:\<G\>\to\<G\>$ be the simplicial homomorphism that is
 the identity on $\<G\>_\0$ and
 zero on $\<1\>$.
 We have the simplicial homomorphisms
 $f\:\<G\>\to G^+$ with $f\circ j=q$ and
 $g\:G^+\to\<G\>$ with $g\circ q\circ i=d\circ j\circ i$.
 We have $f\circ g=\id$.
 Put $D=\<G\>_\0^2\subset\<G\>$.
 Let $k\:D\to\<G\>$ be the inclusion.
 We have the simplicial homomorphism $l\:\<G\>\to D$ such that
 $k\circ l+g\circ f=d$.
 We have $l\circ k=\id$.
 $$
 \xymatrix {
 E
 \ar[r]^-{i} &
 G
 \ar[d]_-{j}
 \ar[dr]^-{q} &
 \\
 D
 \ar[r]^-{k} &
 \<G\>
 \ar@/^2ex/[l]^-{l}
 \ar[r]^-{f} &
 G^+
 \ar@/^2ex/[l]^-{g}
 }
 $$
 The simplicial abelian group $D$ is free.
 By the Freudenthal theorem,
 the map $i\:E\to G$ is $(2n-1)$-connected.
 Since $m\le2n-1$,
 it is $m$-connected.
 Using the Dold--Thom theorem, we see that
 the simplicial homomorphism $\<i\>\:\<E\>\to\<G\>$ is
 $m$-connected.
 One easily sees that
 $(\<i\>,k)\:\<E\>\oplus D\to\<G\>$ is an isomorphism.
 Thus $D$ is $m$-connected.

 For $s\in\N$,
 let $D^{(s)}\subset D$ be the simplicial subgroup equal to:
 $\<G\>_\0^s$ for $s\ge2$, and
 $D$ otherwise.

 \subhead {Decomposition of $D$.}
 Let $r\in\N$, $r\ge2$, be a number.
 By \S~10,
 we have the decomposition $D=D^2\oplus\dotso\oplus D^r$
 where $\<G\>_\0^s=D^s\oplus\dotso\oplus D^r$, $s=2,\dotsc,r$.
 (We have
 $D^s\cong\<E^{\wedge s}\>_\0$ for $s<r$ and
 $D^r=\<G\>_\0^r$.)
 Since $D$ is free and $m$-connected,
 the groups $D^s$ are free and $m$-connected.

 \subhead {The partition $h$.}
 By (14.1),
 for each $s=2,\dotsc,r$,
 there is a partition $(h^s_z\:D^s(\=z)\to D^s(L))_{z\in L}$.
 Combining them,
 we get the partition $(h_z\:D(\=z)\to D(L))_{z\in L}$.
 We have $h_z(D^{(s)}(\=z))\subset D^{(s)}(L)$, $s\in\N$,
 $s\le r$.

 \subhead {The simplicial homomorphism $X$.}
 Let
 $\~G$ be the dummy of $G$,
 $p\:\~G\to G$ be the projection.
 By (13.2),
 $p$ is surjective.
 Thus,
 for the simplicial homomorphism $\<p\>\:\<\~G\>\to\<G\>$,
 we have $\<p\>(\<\~G\>_\0^s)=\<G\>_\0^s$, $s\in\N$.
 Applying (11.1) to each component $D^s$ of the decomposition
 of $D$,
 we get the simplicial homomorphism $X\:D\to\<\~G\>$
 with the following properties:
 \begin {itemize}
 \item [(1)]
 the diagram
 $$
 \xymatrix {
 &&
 \<\~G\>
 \ar[d]^-{\<p\>} \\
 D_{(m)}
 \ar[rru]^-{X|D_{(m)}}
 \ar[rr]^-{\text{\rm inclusion}} &&
 \<G\>
 }
 $$
 is commutative;
 \item [(2)]
 $X(D^{(s)})\subset\<\~G\>_\0^s$, $s\in\N$, $s\le r$.
 \end {itemize}
 We have $\im X\subset\<\~G\>_\0$.

 \subhead {The homomorphism $V$.}
 Let
 $J\:\<G(L)\>\to\<G\>(L)$ be the fusion,
 $R\:\<\~G\>(L)\to\<\~G(L)\>$ be the realization.
 We have the composition
 $$
 \xymatrix {
 V\:D(L)
 \ar[r]^-{X_\#} &
 \<\~G\>(L)
 \ar[r]^-{R} &
 \<\~G(L)\>
 \ar[r]^-{\<p_\#\>} &
 \<G(L)\>.
 }
 $$
 We have $\im V\subset\<\~G(L)\>_\0$.

 \begin {claim} [(15.1)]
 The diagram
 $$
 \xymatrix {
 &
 \<G(L)\>
 \ar[d]^-{J} \\
 D(L)
 \ar[ru]^-{V}
 \ar[r]^-{k_\#} &
 \<G\>(L)
 }
 $$
 is commutative.
 \end {claim}

 \begin {demo} [Proof.]
 Let $\~J\:\<\~G(L)\>\to\<\~G\>(L)$ be the fusion.
 The diagram
 $$
 \xymatrix {
 \<\~G\>(L)
 \ar[r]^-{R} &
 \<\~G(L)\>
 \ar[d]^-{\~J}
 \ar[r]^-{\<p_\#\>} &
 \<G(L)\>
 \ar[d]^-{J} \\
 D(L)
 \ar[u]^-{X_\#}
 \ar[r]^-{X_\#} &
 \<\~G\>(L)
 \ar[r]^-{\<p\>_\#} &
 \<G\>(L)
 }
 $$
 is commutative
 (we invoke (13.4)
 taking into account that $\im X\subset\<\~G\>_\0$).
 We have $J\circ V=\<p\>_\#\circ X_\#=k_\#$
 by the property (1) of $X$.
 \qed
 \end {demo}

 \begin {claim} [(15.2)]
 For $w\in D(L)$,
 we have $\SI(V(w))\subset O_L(\Si(w),1)$.
 \end {claim}

 \begin {demo}
 This follows from (13.5).
 \qed
 \end {demo}

 \begin {claim} [(15.3)]
 We have $V(D^{(s)}(L))\subset\<G(L)\>_\0^s$, $s\in\N$, $s\le r$.
 \end {claim}

 \begin {demo}
 This follows from
 the property (2) of $X$ and
 the claims (13.6) and (9.1).
 \qed
 \end {demo}

 \subhead {The maps $P_z$, $P$.}
 For $z\in L$,
 we have the map $P_z\:G(K)\to\<G(L)\>$,
 $P_z(u)=(V\circ h_z)(l\circ j\circ u\circ e|_{\=z})$.
 We have $P_z(u)\in\<G(L)\>_\0$
 since $\im V\subset\<\~G(L)\>_\0$.
 We have $\SI(P_z(u))\subset O_L(z,b_1)$
 (by
 the definition of a partition,
 the claim (15.2), and
 the inequality $b_1\ge2$).

 We have the map $P\:G(K)\to\<G(L)\>$,
 $P(u)=V(l\circ j\circ u\circ e)$.
 We have
 $$
 \sum_{z\in L}P_z(u)=P(u).
 $$

 \subhead {The homomorphism $M$.}
 We have the additive homomorphism $M\:\<G(K)\>\to\<G(L)\>$,
 $$
 M(\`u\')=
 \sum_{Z\subset L\%\Ep_L(Z)\ge b_3}
 (-1)^{\#Z}
 \`u\circ e_Z\'
 \prod_{z\in Z}
 P_z(u).
 $$
 Here and in all our $\prod$'s,
 we mean that
 the order of factors is induced by some fixed order on $L$.
 (Moreover, one can see that
 the factors commute everywhere.)

 \begin {claim} [(15.4)]
 For $U\in\<G(K)\>$,
 we have $\Th(M(U))\ge\min(\Th(U)+1,\Et(U))$.
 \end {claim}

 \begin {demo} [Proof.]
 Suppose that
 $\Th(U)\ge s-1$ and
 $\Et(U)\ge s$ ($s\in\N_+$).
 We show that
 $\Th(M(U))\ge s$.
 Take a set $T\subset L$ with $\#T<s$.
 We show that
 $M(U)\|_T=0$.

 {\it The case $\Ep_L(T)\ge b_4$.\/}
 Put $I=\{Z\subset L:\Ep_L(Z)\ge b_3\}$.
 For $u\in G(K)$,
 we have
 $$
 M(\`u\')\|_T=
 \sum_{Z\in I}
 (-1)^{\#Z}
 \`u\circ e_Z\'\|_T
 \prod_{z\in Z}
 P_z(u)\|_T.
 $$
 The sets $O_L(y,b_1)$, $y\in T$, ({\it balls}) do not
 intersect.
 Moreover, the distance ($\Rh_L$) between simplices of distinct
 balls is at least $b_3$
 (since $b_4\ge2b_1+b_3$).
 The distance between simplices of a ball is smaller than $b_3$
 (since $b_3\ge2b_1$).
 Let $I_0$ be the set of sets $Z\subset L$ that
 are contained in the union of the balls and
 have at most one simplex in each ball.
 Show that
 our sum over $Z\in I$ equals the same sum but over $Z\in I_0$.
 We have $I_0\subset I$.
 If $Z\in I\setminus I_0$,
 there is a simplex $z\in Z\setminus O_L(T,b_1)$;
 then $P_z(u)\|_T=0$
 because:
 $P_z(u)\in\<G(L)\>_\0$,
 $\SI(P_z(u))\subset O_L(z,b_1)$, and
 $O_L(z, b_1)\cap T=\varnothing$.
 Thus the corresponding summand is zero.

 Put
 $$
 I_0'=\coprod_{S\subset T}W_S,
 $$
 where $W_S$ is the set of maps $w\:S\to L$ such that
 $w(y)\in O_L(y,b_1)$, $y\in T$.
 We have the bijection $I_0'\to I_0$, $(S,w)\mapsto w(S)$.
 Thus
 $$
 M(\`u\')\|_T=
 \sum_{(S,w)\in I_0'}
 (-1)^{\#S}
 \`u\circ e_{w(S)}\'\|_T
 \prod_{y\in S}
 P_{w(y)}(u)\|_T.
 $$

 For $y\in T$,
 let $t_y\:G(\=y)\to G_T$ be the canonical monomorphism of a
 factor to a product.
 Show that
 for $(S,w)\in I_0'$,
 $$
 (u\circ e_{w(S)})\|_T=
 \prod_{y\in T\setminus S}
 t_y(u\circ e|_{\=y}).
 $$
 If $y\in S$,
 we have
 $y\in O_L(w(y),b_1)$, and
 $e_{w(S)}$ sends the simplex $y$ to a vertex of $K$;
 then $u\circ e_{w(S)}|_{\=y}=1$
 since $G_0=1$.
 If $y\in T\setminus S$,
 we have
 $y\notin O_L(w(S),b_2)$
 (since $b_4\ge b_1+b_2$), and
 $e_{w(S)}|_{\=y}=e|_{\=y}$.
 Thus we have the desired equality.

 For $(S,w)\in I_0'$ and $y\in S$,
 we have $P_{w(y)}(u)\|_T=\<t_y\>(P_{w(y)}(u)|_{\=y})$.
 This is because
 $\SI(P_{w(y)}(u))\subset O_L(w(y),b_1)$ and
 $O_L(w(y),b_1)\cap T=\{y\}$
 (since $b_4\ge2b_1$).

 Thus
 \begin {multline*}
 M(\`u\')\|_T=
 \sum_{(S,w)\in I_0'}
 (-1)^{\#S}
 \bigl(
 \prod_{y\in T\setminus S}
 \`t_y(u\circ e|_{\=y})\'
 \bigr)
 \bigl(
 \prod_{y\in S}
 \<t_y\>(P_{w(y)}(u)|_{\=y})
 \bigr)= \\ =
 \sum_{S\subset T}
 (-1)^{\#S}
 \bigl(
 \prod_{y\in T\setminus S}
 \`t_y(u\circ e|_{\=y})\'
 \bigr)
 \bigl(
 \sum_{w\in W_S}
 \prod_{y\in S}
 \<t_y\>(P_{w(y)}(u)|_{\=y})
 \bigr)= \\ =
 \sum_{S\subset T}
 (-1)^{\#S}
 \bigl(
 \prod_{y\in T\setminus S}
 \<t_y\>(\`u\circ e\'|_{\=y})
 \bigr)
 \bigl(
 \prod_{y\in S}
 \sum_{z\in O_L(y,b_1)}
 \<t_y\>(P_z(u)|_{\=y})
 \bigr)= \\ =
 \prod_{y\in T}
 \<t_y\>\bigl(
 \`u\circ e\'|_{\=y}-
 \sum_{z\in O_L(y,b_1)}
 P_z(u)|_{\=y}
 \bigr).
 \end {multline*}
 We may extend the domain of the last sum to $z\in L$ because
 for $z\in L\setminus O_L(y,b_1)$,
 we have $P_z(u)|_{\=y}=0$
 because:
 $P_z(u)\in\<G(L)\>_\0$,
 $\SI(P_z(u))\subset O_L(z,b_1)$, and
 $O_L(z,b_1)\cap\=y=\varnothing$ for such $z$.
 We have
 $$
 M(\`u\')\|_T=
 \prod_{y\in T}
 \<t_y\>(\`u\circ e\'|_{\=y}-P(u)|_{\=y}).
 $$

 For $y\in T$,
 let $J_y\:\<G(\=y)\>\to\<G\>(\=y)$ be the fusion.
 Obviously, it is an isomorphism.
 We have the commutative diagram
 $$
 \xymatrix {
 &
 \<G(L)\>
 \ar[d]^-{J}
 \ar[r]^-{\?|_{\=y}} &
 \<G(\=y)\>
 \ar[d]^-{J_y} \\
 D(L)
 \ar[ur]^-{V}
 \ar[r]^-{k_\#} &
 \<G\>(L)
 \ar[r]^-{\?|_{\=y}} &
 \<G\>(\=y)
 }
 $$
 (we invoke (15.1)).
 We have $J_y(\`u\circ e\'|_{\=y}-P(u)|_{\=y})=
 J_y(\`u\circ e\'|_{\=y}-V(l\circ j\circ u\circ e)|_{\=y})=
 j\circ u\circ e|_{\=y}-k\circ l\circ j\circ u\circ e|_{\=y}=
 1+g\circ f\circ j\circ u\circ e|_{\=y}=
 1+g\circ q\circ u\circ e|_{\=y}$.
 We have the homomorphism $a_y\:G_K\to\<G_T\>$ (in the additive
 group), $a_y(v)=
 (\<t_y\>\circ J_y^{-1})((g\circ q\circ v\circ e)_y)$.
 We have
 $$
 M(\`u\')\|_T=
 \prod_{y\in T}
 (1+a_y(u\|_K)).
 $$
 Since $\Et(U)>\#T$,
 by (8.1),
 $M(U)\|_T=0$.

 {\it The converse case.\/}
 There are distinct simplices $y_0,y_1\in T$ with
 $\Rh_L(y_0,y_1)<b_4$.
 For each $y\in T\setminus\{y_1\}$,
 consider the simplex $x\in K$ spanned by the set
 $e(O_L(y,b_5))$.
 Let $S\subset K$ be the set of these simplices.
 We have $\#S<s-1$.
 For each $y\in T$,
 there exists a simplex $y'\in T\setminus\{y_1\}$ such that
 $O_L(y,2b_2)\subset O_L(y',b_5)$:
 we may let $y'$ be equal to:
 $y_0$ if $y=y_1$, and
 $y$ otherwise
 (we use the inequality $b_5\ge2b_2+b_4$).
 Thus,
 for every $y\in T$,
 there exists a simplex $x\in S$ such that
 $e(O_L(y,2b_2))\subset\=x$.
 Let $e'\:\ovl{O_L(T,b_1)}\to\-S$ be the abridgement of $e$
 (we use the inequality $b_1\le2b_2$).

 Take a set $Z\subset L$ such that $\Ep_L(Z)\ge b_3$.
 Show that $e_Z(\-T)\subset\-S$.
 It suffices to check that
 $e_Z(y)\in\-S$ for $y\in T$.
 If $y\notin O_L(Z,b_2)$,
 then $e_Z(y)=e(y)\in\-S$.
 Otherwise,
 $y\in O_L(z,b_2)$ for some $z\in Z$.
 Then $e_Z(y)=e_z(y)\in\=x$, where
 $x\in K$ is the simplex spanned by $e(O_L(z,b_2))$.
 We have $e(O_L(z,b_2))\subset e(O_L(y,2b_2))$
 Thus $e_Z(y)\in\-S$.
 Let $\~e_Z\:\-T\to\-S$ be the abridgement of $e_Z$.

 We have the additive homomorphism $\~M\:\<G(\-S)\>\to\<G(\-T)\>$,
 $$
 \~M(\`\~u\')=
 \sum_{Z\subset O_L(T,b_1)\%\Ep_L(Z)\ge b_3}
 (-1)^{\#Z}
 \`\~u\circ\~e_Z\'
 \prod_{z\in Z}
 (V\circ h_z)(l\circ j\circ\~u\circ e'|_{\=z})|_{\-T}.
 $$
 Show that the diagram
 $$
 \xymatrix {
 \<G(K)\>
 \ar[d]_-{\?|_{\-S}}
 \ar[r]^-{M} &
 \<G(L)\>
 \ar[d]_-{\?|_{\-T}} \\
 \<G(\-S)\>
 \ar[r]^-{\~M} &
 \<G(\-T)\>
 }
 $$
 is commutative.
 We have
 $$
 M(\`u\')|_{\-T}=
 \sum_{Z\subset L\%\Ep_L(Z)\ge b_3}
 (-1)^{\#Z}
 \`u\circ e_Z\'|_{\-T}
 \prod_{z\in Z}
 P_z(u)|_{\-T}.
 $$
 The summands with $Z\not\subset O_L(T,b_1)$ equal zero
 (if $z\in Z\setminus O_L(T,b_1)$,
 then $P_z(u)|_{\-T}=0$
 because:
 $P_z(u)\in\<G(L)\>_\0$,
 $\SI(P_z(u))\subset O_L(z,b_1)$, and
 $O_L(z, b_1)\cap\-T=\varnothing$).
 We get
 \begin {multline*}
 M(\`u\')|_{\-T}=
 \sum_{Z\subset O_L(T,b_1)\%\Ep_L(Z)\ge b_3}
 (-1)^{\#Z}
 \`u\circ e_Z\'|_{\-T}
 \prod_{z\in Z}
 (V\circ h_z)(l\circ j\circ u\circ e|_{\=z})|_{\-T}= \\ =
 \~M(\`u\'|_{\-S}).
 \end {multline*}

 Since $\Th(U)>\#S$,
 $U\|_S=0$.
 Thus
 $U|_{\-S}=0$.
 We get $M(U)|_{\-T}=\~M(U|_{\-S})=0$.
 Thus
 $M(U)\|_T=0$.
 \qed
 \end {demo}

 \begin {claim} [(15.5)]
 For $U\in\<G(K)\>$,
 we have $\Et(M(U))\ge\min(\Et(U),r)$.
 \end {claim}

 \begin {demo} [Proof.]
 We have the additive homomorphism $N\:\<G_K\>\to\<G_L\>$,
 $$
 N(\`v\')=
 \sum_{Z\subset L\%\Ep_L(Z)\ge b_3}
 (-1)^{\#Z}
 \`v\circ e_Z\'
 \prod_{z\in Z}
 (V\circ h_z)((l\circ j\circ v\circ e)_z)\|_L.
 $$
 The diagram
 $$
 \xymatrix {
 \<G(K)\>
 \ar[d]_-{\?\|_K}
 \ar[r]^-{M} &
 \<G(L)\>
 \ar[d]_-{\?\|_L} \\
 \<G_K\>
 \ar[r]^-{N} &
 \<G_L\>
 }
 $$
 is commutative.
 It suffices to show that
 $N$ is $r$-strict.
 For $z\in L$,
 we have
 the homomorphism $t_z\:G_K\to G(\=z)$, $t_z(v)=(v\circ e)_z$,
 and
 the additive homomorphism $B_z\:\<G(\=z)\>\to\<G(L)\>$,
 $B_z(\`v\')=(V\circ h_z)(l\circ j\circ v)$.
 We have the homomorphisms
 $e_Z^\#\:G_K\to G_L$ and
 $\?\|_L\:G(L)\to G_L$.
 We have
 (for $v\in G_K$)
 $$
 N(\`v\')=
 \sum_{Z\subset L\%\Ep_L(Z)\ge b_3}
 (-1)^{\#Z}
 \<e_Z^\#\>(\`v\')
 \prod_{z\in Z}
 (B_z\circ\<t_z\>)(\`v\')\|_L.
 $$
 By (9.1) and (9.2),
 it suffices to show that
 the homomorphisms $B_z$ are $r$-strict.
 The homomorphism $B_z$ equals the composition
 $$
 \xymatrix {
 \<G(\=z)\>
 \ar[r]^-{J_z} &
 \<G\>(\=z)
 \ar[r]^-{l_\#} &
 D(\=z)
 \ar[r]^-{h_z} &
 D(L)
 \ar[r]^-{V} &
 \<G(L)\>,
 }
 $$
 where $J_z$ is the fusion.
 For $s\in\N$,
 we have:
 $J_z(\<G(\=z)\>_\0^s)=\<G\>_\0^s(\=z)$;
 $l_\#(\<G\>_\0^s(\=z))=D^{(s)}(\=z)$
 (since $l$ is identical on $D$);
 $h_z(D^{(s)}(\=z))\subset D^{(s)}(L)$ for $s\le r$
 (a property of the partition $h$);
 $V(D^{(s)}(L))\subset\<G(L)\>_\0^s$ for $s\le r$
 (by (15.3)).
 Thus
 $B_z(\<G(\=z)\>_\0^s)\subset\<G(L)\>_\0^s$ for $s\le r$,
 which is what we need.
 \qed
 \end {demo}

 Put $Q=|K|$ ($=|L|$).

 \begin {claim} [(15.6)]
 For $U\in\<G(K)\>$,
 we have $[|M(U)|]=[|U|]$ in the ring $\<[Q,|G|]\>$.
 \end {claim}

 \begin {demo} [Proof.]
 Take $u\in G(K)$ and $z\in L$.
 We have $P_z(u)\in\<G(L)\>_\0$.
 By the construction of $P_z$,
 all the sections in the ensemble $P_z(u)$
 lift to $\~G$.
 By (13.1),
 the space $|\~G|$ is contractible.
 Thus $[|P_z(u)|]=0$.
 Applying the ring homomorphism
 $[|\?|]\:\<G(L)\>\to\<[Q,|G|]\>$ to the equality defining $M$,
 we get $[|M(\`u\')|]=[|\`u\circ e\'|]=[|\`u\'|]$
 since $|e|\:Q\to Q$ is homotopic to the identity.
 \qed
 \end {demo}


 \head {16. Main procedure}

 Let
 $K$ be a polyhedron with $\dim K\le m$ ($m\in\N$) and
 $E$ be an $(n-1)$-connected ($n\in\N$) simplicial set with a
 single vertex.
 Suppose that $m\le2n-1$.
 Put
 $Q=|K|$ and
 $G=FE$.

 \begin {claim} [(16.1)]
 Let $U\in\<G(K)\>$ be an ensemble with $\Et(U)\ge s$
 ($s\in\N$).
 Then there exist
 a polyhedron $L$ with the body $Q$ and
 an ensemble $V\in\<G(L)\>$ with
 $\Th(V)\ge s$ and
 $[|V|]=[|U|]$ in $\<[Q,|G|]\>$.
 \end {claim}

 \begin {demo} [Proof.]
 To get the desired pair $(L,V)$,
 take the pair $(K,U)$ and
 apply the pair $(\DE^c,M)$ of operations of \S~15 $s$ times.
 We put $r=s$.
 The desired properties follow from (15.4), (15.5), and (15.6).
 \qed
 \end {demo}


 \head {17. The function $\Th$: the topological version}

 Let $X$ and $Y$ be spaces.
 For $A\in\<C(X,Y)\>$,
 put $\Th(A)=
 \inf\{\#V:\text{\rm finite}\ V\subset X,\ A|_V\ne0\}\in\^\N$.

 Let
 $X'$ and $Y'$ be spaces,
 $g\:X'\to X$ and $h\:Y\to Y'$ be continuous maps.
 We have the map $t\:C(X,Y)\to C(X',Y')$,
 $t(a)=h\circ a\circ g$.
 We have the homomorphism $\<t\>\:\<C(X,Y)\>\to\<C(X',Y')\>$.

 \begin {claim} [(17.1)]
 For $A\in\<C(X,Y)\>$,
 we have $\Th(\<t\>(A)\>\ge\Th(A)$.
 \end {claim}

 \begin {demo} [Proof.]
 Take a finite $V'\subset X'$ with $\#V'<\Th(A)$.
 We show that $\<t\>(A)|_{V'}=0$.
 Put $V=g(V')\subset X$.
 We have $\#V<\Th(A)$.
 Thus $A|_V=0$.
 Let $\~g\:V'\to V$ be the abridgement of $g$.
 We have the map $\~t\:C(V,Y)\to C(V',Y')$,
 $\~t(\~a)=h\circ\~a\circ\~g$.
 The diagram
 $$
 \xymatrix {
 C(X,Y)
 \ar[d]_-{\?|_V}
 \ar[r]^-{t} &
 C(X',Y')
 \ar[d]^-{\?|_{V'}} \\
 C(V,Y)
 \ar[r]^-{\~t} &
 C(V',Y')
 }
 $$
 is commutative.
 We have $\<t\>(A)|_{V'}=\<\~t\>(A|_V)=0$.
 \qed
 \end {demo}

 \subhead {A characterization of the order.}
 Let
 $U$ be an abelian group and
 $f\:[X,Y]\to U$ be a map.
 We have the homomorphism $\-f\:\<[X,Y]\>\to U$,
 $\-f(\`w\')=f(w)$.

 \begin {claim} [(17.2)]
 The condition $\ord f\le r$ ($r\in\N$) is equivalent to the
 condition that $\-f([A])=0$ for every $A\in\<C(X,Y)\>$ with
 $\Th(A)>r$.
 \end {claim}

 \begin {demo} [Proof.]
 Let $E_r$, $I_r$, and $D_r$ be as in \S~1.
 We have the homomorphism $h\:\<C(X,Y)\>\to D_r$,
 $h(\`a\')=I_r(a)$.
 It is surjective.
 One easily sees that
 for $A\in\<C(X,Y)\>$,
 the conditions $h(A)=0$ and $\Th(A)>r$ are equivalent.
 We have the homomorphism $\~f\:\<C(X,Y)\>\to U$,
 $\~f(A)=\-f([A])$.
 The condition $\ord f\le r$ is equivalent to the existence of
 a homomorphism $l\:D_r\to U$ with $l\circ h=\~f$.
 The latter is equivalent to the condition $\~f|\ker h=0$,
 that is, the condition that $\-f([A])=0$ for every
 $A\in\<C(X,Y)\>$ with $\Th(A)>r$.
 \qed
 \end {demo}


 \head {18. Geometric realization and simplicial approximation}

 Let
 $K$ be a polyhedron and
 $E$ be a simplicial set.
 Put $Q=|K|$.

 \begin {claim} [(18.1)]
 For $U\in\<E(K)\>$,
 we have $\Th(|U|)=\Th(U)$.
 \qed
 \end {claim}

 \begin {claim} [(18.2)]
 Let $B\in\<C(Q,|E|)\>$ be an ensemble.
 Then there exist
 a polyhedron $L$ with the body $Q$ and
 an ensemble $V\in\<E(L)\>$ with
 $\Th(V)\ge\Th(B)$ and
 $[|V|]=[B]$ in $\<Q,|E|\>$.
 \end {claim}

 \begin {demo} [Proof.]
 There are
 a finite set $I$,
 a map $k\:I\to C(Q,|E|)$, and
 an element $g\in\<I\>$
 such that $\<k\>(g)=B$.
 Put $b_i=k(i)$, $i\in I$.
 For $q\in Q$,
 we have the equivalence $R_q=\{(i,j):b_i(q)=b_j(q)\}$ on $I$.
 For a finite set $W\subset Q$,
 put
 $$
 R_W=\bigcap_{q\in W}R_q.
 $$
 The map $i\mapsto b_i|_W$ is subordinate to the equivalence
 $R_W$
 (that is, constant on the classes of $R_W$).
 We have the commutative diagram
 $$
 \xymatrix {
 I
 \ar[d]_-{p_W}
 \ar[r]^-{k} &
 C(Q,|E|)
 \ar[d]^-{\?|_W} \\
 I/R_W
 \ar[r]^-{k_W} &
 C(W,|E|),
 }
 $$
 where $p_W$ is the projection.
 The map $k_W$ is injective.
 We have $\<k_W\>(\<p_W\>(g))=\<k\>(g)|_W=B|_W$.
 If $\#W<\Th(B)$,
 then $B|_W=0$, and
 thus $\<p_W\>(g)=0$.

 We have the continuous map $b=(b_i)_{i\in I}\:Q\to|E|^I$.
 Let $h\:|E^I|\to|E|^I$ be the canonical continuous bijection.
 Since
 $I$ is finite and
 $Q$ is Haudorff and compact,
 the map $c=h^{-1}\circ b\:Q\to|E^I|$ is continuous.

 To each equivalence $R$ on $I$ assign the simplicial subset
 $D(R)\subset E^I$,
 $D(R)_n=
 \{(e_i)_{i\in I}\in E_n^I:
 (i,j)\in R\Rightarrow
 e_i=e_j\}$ (the diagonal).
 For $q\in Q$,
 we have $c(q)\in|D(R_q)|\subset|E^I|$.
 We have the simplicial subset $M\subset E^I$,
 $$
 M=\bigcup_{q\in Q}D(R_q).
 $$
 We have $c(Q)\subset|M|\subset|E^I|$.
 Let $c'\:Q\to|M|$ be the abridgement of $c$.
 By the simplicial approximation theorem,
 there are
 a polyhedron $L$ with the body $Q$ and
 a section $u'\in M(L)$ such that
 the map $|u'|\:Q\to|M|$ is homotopic to $c'$.
 Let $u\in E^I(L)$ be the composition of $u'$ and the inclusion
 $M\to E^I$.
 We have $u=(u_i)_{i\in I}$, where $u_i\in E(L)$.
 The map $|u_i|\:Q\to|E|$ is homotopic to $b_i$.
 We have the map $l\:I\to E(L)$, $l(i)=u_i$.
 Put $V=\<l\>(g)$.
 We have $[|V|]=[B]$.

 For a simplex $y\in L$, $\dim y=s$,
 we have $u_s(y)\in M_s$,
 that is, there is a point $q=q_y\in Q$ such that
 $u_s(y)\in D(R_q)_s$,
 that is, $u_i|_{\=y}=u_j|_{\=y}$ for $(i,j)\in R_q$,
 that is, the map $i\mapsto u_i|_{\=y}$ is subordinate to
 $R_q$.

 Take a set $T\subset L$.
 Put $W=\{q_y:y\in T\}$.
 We have $\#W\le\#T$.
 The map $i\mapsto u_i\|_T$ is subordinate to $R_W$.
 We have the commutative diagram
 $$
 \xymatrix {
 I
 \ar[d]_-{p_W}
 \ar[r]^-{l} &
 E(L)
 \ar[d]^-{\?\|_T} \\
 I/R_W
 \ar[r]^-{l_T} &
 E_T.
 }
 $$
 We have $V\|_T=\<l\>(g)\|_T=\<l_T\>(\<p_W\>(g))$.
 If $\#T<\Th(B)$,
 then:
 $\#W<\Th(B)$,
 $\<p_W\>(g)=0$, and
 $V\|_T=0$.
 Thus $\Th(V)\ge\Th(B)$.
 \qed
 \end {demo}


 \head {19. Some subgroups of $\<[Q,|G|]\>$.}

 Let
 $Q$ be a polyhedral body, $\dim Q\le m$ ($m\in\N$), and
 $E$ be a $(n-1)$-connected ($n\in\N$) simplicial set with a
 single vertex.
 Suppose that $m\le2n-1$.
 Put $G=FE$.
 Define the subgroups $P,M_s,J_s\subset\<C(Q,|G|)\>$, $s\in\N$:
 put
 $P=\<C(Q,|G_{(m)}|)\>$
 (we have $C(Q,|G_{(m)}|)\subset C(Q,|G|)$),
 $M_s=\{B:\Th(B)\ge s\}$, and
 let $J_s$ be generated by all elements of the form
 $(\`b_1\'-1)\dotso(\`b_k\'-1)$, where
 $k\in\N$,
 $b_l\in C(Q,|\Ga_{s_l}G|)\subset C(Q,|G|)$, and
 $s_1+\dotso+s_k\ge s$.
 ($M_s$ and $J_s$ are ideals.
 Conjecture: $M_s\subset J_s$.)
 For a subgroup $S\subset\<C(Q,|G|)\>$,
 let $[S]\subset\<[Q,|G|]\>$ be its image under the
 homomorphism
 $[\?]\:\<C(Q,|G|)\>\to\<[Q,|G|]\>$.

 \begin {claim} [(19.1)]
 For $s\in\N$,
 we have $[M_s]=[P\cap M_s]=[J_s]$.
 \end {claim}

 \begin {demo} [Proof.]
 {\it The inclusion $[M_s]\subset[J_s]$.\/}
 Take an element $B\in M_s$.
 We have $\Th(B)\ge s$.
 By (18.2),
 there are
 a polyhedron $L$ with the body $Q$ and
 an ensemble $V\in\<G(L)\>$ with
 $\Th(V)\ge s$ and
 $[|V|]=[B]$.
 It suffices to show that
 $|V|\in J_s$.
 By (7.2),
 $\Et(V)\ge s$.
 Let $I_s\subset\<G(L)\>$ be, as in \S~7, the subgroup
 generated by all elements of the form
 $(\`v_1\'-1)\dotso(\`v_k\'-1)$, where
 $k\in\N$,
 $v_l\in(\Ga_{s_l}G)(L)\subset G(L)$, and
 $s_1+\dotso+s_k\ge s$.
 By (7.1),
 $V\in I_s$.
 Obviously,
 $|V|\in J_s$.

 {\it The inclusion $[P\cap M_s]\supset[J_s]$.\/}
 Take an element $B\in\<C(Q,|G|)\>$,
 $B=(\`b_1\'-1)\dotso(\`b_k\'-1)$, where
 $k\in\N$,
 $b_l\in C(Q,|\Ga_{s_l}G|)\subset C(Q,|G|)$, and
 $s_1+\dotso+s_k\ge s$.
 Such elements generate $J_s$.
 Thus it suffices to show that
 $[B]\in[P\cap M_s]$.
 Choose a polyhedron $K$ with the body $Q$.
 Since $\Ga_sG$ are Kan sets,
 there are sections $u_l\in(\Ga_{s_l}G)(K)$ with
 $[|u_l|]=[b_l]$ in $[Q,|G|]$.
 Put $U=(\`u_1\'-1)\dotso(\`u_k\'-1)\in\<G(L)\>$.
 We have $[|U|]=[B]$ in $\<[Q,|G|]\>$.
 By (7.1),
 $\Et(U)\ge s$.
 By (16.1),
 there are
 a polyhedron $L$ with the body $Q$ and
 an ensemble $V\in\<G(L)\>$ with
 $\Th(V)\ge s$ and
 $[|V|]=[|U|]$ in $\<[Q,|G|]\>$.
 Obviously,
 $|V|\in P$.
 By (17.1),
 $\Th(|V|)\ge s$.
 Thus,
 $[B]=[|V|]$ and
 $|V|\in P\cap M_s$.
 \qed
 \end {demo}


 \head {20. Step from $[Q,|G|]$ to $[X,Y]$}

 Let
 $X$ be a finite CW-complex, $\dim X\le m$ ($m\in\N$), and
 $Y$ be an $(n-1)$-connected ($n\in\N$) CW-complex.
 Suppose that $m<2n-1$.
 We have the subgroups $L_s\subset\<C(X,Y)\>$, $s\in\N$:
 $L_s=\{A\:\Th(A)\ge s\}$.
 Let $B=(B_s)_{s=1}^\infty$ be the Curtis filtration of
 $[X,Y]$.
 For $s\in\N$,
 we have the subgroup $H_s\subset\<[X,Y]\>$ generated by all
 elements of the form $(\`w_1\'-1)\dotso(\`w_k\'-1)$, where
 $k\in\N$,
 $w_l\in B_{s_l}$, and
 $s_1+\dotso+s_k\ge s$.
 (It is an ideal.)
 For a subgroup $R\subset\<C(X,Y)\>$,
 let $[R]\subset\<[X,Y]\>$ be its image under the homomorphism
 $[\?]\:\<C(X,Y)\>\to\<[X,Y]\>$.

 \begin {claim} [(20.1)]
 We have $[L_s]=H_s$,
 $s\in\N$.
 \end {claim}

 \begin {demo} [Proof.]
 There are
 a polyhedral body $Q$, $\dim Q\le m$, and
 a homotopy equivalence $g\:Q\to X$.
 Let $g'\:X\to Q$ be a homotopy inverse map.
 There are
 a simplicial set $E$ with a single vertex and
 a homotopy equivalence $k\:Y\to|E|$.
 Put $G=FE$.
 Let $i\:E\to G$ be the canonical simplicial map.
 By the Freudenthal theorem,
 it is $(2n-1)$-connected.
 The map $h=|i|\circ k\:Y\to|G|$ is also $(2n-1)$-connected.
 Since $m\le 2n-1$,
 there is a map $h'\:|G_{(m)}|\to Y$ such that
 the map $h\circ h'$ is homotopic to the inclusion
 $|G_{(m)}|\to|G|$.
 We have the map $t\:C(X,Y)\to C(Q,|G|)$,
 $t(a)=h\circ a\circ g$.
 Since $m<2n-1$,
 it induces an isomorphism $\-t\:[X,Y]\to[Q,|G|]$.
 We have the map $t'\:C(Q,|G_{(m)}|)\to C(X,Y)$,
 $t'(b)=h'\circ b\circ g'$.
 For $b\in C(Q,|G_{(m)}|)\subset C(Q,|G|)$,
 we have $[t'(b)]=\-t^{-1}([b])$.
 One can see that
 \begin {equation}
 \-t(B_s)=
 \{
 [b]\in[Q,|G|]:
 b\in C(Q,|\Ga_sG|)\subset C(Q,|G|)
 \},
 \qquad s\in\N.
 \tag {$*$}
 \end {equation}

 Let $P,M_s,J_s\subset\<C(Q,|G|)\>$ be as in \S~19.
 We have the homomorphisms
 $\<t\>\:\<C(X,Y)\>\to\<C(Q,|G|)\>$ and
 $\<t'\>\:P=\<C(Q,|G_{(m)}|)\>\to\<C(X,Y)\>$.
 By (17.1),
 $\<t\>(L_s)\subset M_s$, and
 $\<t'\>(P\cap M_s)\subset L_s$.
 We have the ring isomorphism
 $\<\-t\>\:\<[X,Y]\>\to\<[Q,|G|]\>$.
 It follows from \thetag {$*$} that
 $\<\-t\>(H_s)=[J_s]$.
 Using (19.1),
 we get
 $\<\-t\>([L_s])=[\<t\>(L_s)]\subset[M_s]=[J_s]=\<\-t\>(H_s)$.
 Hence $[L_s]\subset H_s$, and
 $[L_s]\supset[\<t'\>(P\cap M_s)]=\<\-t^{-1}\>([P\cap M_s])=
 \<\-t\>^{-1}([J_s])=H_s$.
 \qed
 \end {demo}

 \begin {demo} [Proof of Theorem (1.1).]
 We have the homomorphism $\-f\:\<[X,Y]\>\to U$,
 $\-f(\`w\')=f(w)$.
 By (17.2),
 the condition $\ord f<s$ ($s\in\N_+$) is equivalent to the
 condition $\-f|[L_s]=0$.
 Obviously,
 the condition $\deg_Bf<s$ is equivalent to the condition
 $\-f|H_s=0$.
 Now note that $[L_s]=H_s$
 by (20.1).
 \qed
 \end {demo}



 {\noindent \tt ssp@pdmi.ras.ru}

 {\noindent \tt http://www.pdmi.ras.ru/\tildeaccent{}ssp}

 \end {document}